\documentclass[12pt]{article}

\usepackage{amsmath}
\usepackage{amsfonts}
\usepackage{amsthm}
\usepackage{hyperref}
\usepackage{enumitem}
\usepackage{array}
\usepackage[nottoc,numbib]{tocbibind}

\newcommand{\dimh}{\mbox{$\dim_{\mathrm{H}}$}}
\newcommand{\dm}{\mbox{$\mu^{\mathrm{H}}$}}

\begin{document}


\author{Attila Losonczi}
\title{The Hausdorff-integral and its applications}

\date{\today}

\newtheorem{thm}{\qquad Theorem}[section]
\newtheorem{prp}[thm]{\qquad Proposition}
\newtheorem{lem}[thm]{\qquad Lemma}
\newtheorem{cor}[thm]{\qquad Corollary}
\newtheorem{rem}[thm]{\qquad Remark}
\newtheorem{ex}[thm]{\qquad Example}
\newtheorem{df}[thm]{\qquad Definition}
\newtheorem{prb}{\qquad Problem}

\newtheorem*{thm2}{\qquad Theorem}
\newtheorem*{cor2}{\qquad Corollary}
\newtheorem*{prp2}{\qquad Proposition}

\maketitle

\begin{abstract}

\noindent

We present a new type of integral that is supposed to extend the usability of the Lebesgue integral in certain types of investigations. It is based on the Hausdorff dimension and measure.
We examine the basic properties of the integral and its similarities to the properties of the Lebesgue integral.
We present many applications as well.

\noindent
\footnotetext{\noindent
AMS (2020) Subject Classifications:  28A25, 28A78\\

}

\end{abstract}


\section{Introduction}

The aim of this paper is to define a new type of integral in order to refine some properties of the Lebesgue integral in certain types of investigations. Our main motivation is the following. In analysis, there are many notions that are based on a kind of measurement of non-negative function values, while some other notions are based on non-zero function values. E.g. several types of (pseudo-)metrics on function spaces, or deficiency measurement (see \cite{lamd}), which is we want to measure how much a certain property is not fulfilled for some objects. For example, we want to measure how much a given function is not continuous on an interval, or how much it is not even/odd, not periodic, etc. If one applies the Lebesgue measure and integral in those cases, then one can get reasonable results, however in many cases, one ends up with $0$ as a result, because the assigned descriptor function is only non-zero on a set of Lebesgue \mbox{measure $0$}. Our new tool can penetrate into those realms too, in other words, it can provide a wide range of measurements in those cases as well. 




\smallskip

The structure of the paper is the following. First we define some algebraic structures on $[0,+\infty)\times[-\infty,+\infty]$ as they will play a central role in the bases of the integral later. Then we define the integral and some related notions, e.g. some generalized metrics as well. We give a new type of set measurement that can also be derived as a Hausdorff integral. As natural, we build the basics of the theory of the new integral, we provide all the necessary tools that will be needed in further investigations. Afterwards, we examine if the analog of the Beppo-Levi, Fatou, Riesz-Fischer theorem holds, and if yes, under which circumstances. At the end, we add many applications, we define some deficiency measurements in concrete cases that are based on our new tools.

\subsection{Basic notions and notations}\label{sbnn}

Here we enumerate the basics that we will apply throughout the paper.

If $A,B$ are sets, then $A\triangle B=(A-B)\cup(B-A)$.

\smallskip

For $K\subset\mathbb{R}$, $\chi_K$ will denote the characteristic function of $K$ i.e. $\chi_K(x)=1$ if $x\in K$, otherwise $\chi_K(x)=0$.

\smallskip

We will use the usual notations $\mathbb{N}_0=\mathbb{N}\cup\{0\},\ \overline{\mathbb{N}}_0=\mathbb{N}\cup\{0,+\infty\},\newline \mathbb{R}^+=\{x\in\mathbb{R}:x>0\},\ \mathbb{R}^+_0=\mathbb{R}^+\cup\{0\},\ \overline{\mathbb{R}}^+_0=\mathbb{R}^+_0\cup\{+\infty\},\ \overline{\mathbb{R}}=\mathbb{R}\cup\{-\infty,+\infty\}.$

\smallskip

Some usual operations and relation with $\pm\infty$: $(+\infty)+(+\infty)=+\infty,\ (-\infty)+(-\infty)=-\infty$; if $r\in\mathbb{R}$, then $r+(+\infty)=+\infty,\ r+(-\infty)=-\infty,\ {-\infty<r<+\infty}$. $+\infty+(-\infty)$ is undefined.

\smallskip


Set $S(x,\delta)=\{y \in\mathbb{R}^n:d(x,y)<\delta\},\ \ \overset{.}{S}(x,\delta)=S(x,\delta)-\{x\}\ \ (x\in\mathbb{R}^n,\delta>0)$.

\smallskip

In $\mathbb{R}^n\ \pi_k$ denotes the \textbf{projection} to the $k^{th}$ coordinate, in other words $\pi_k(x_1,\dots,x_n)=x_k\ \ (1\leq k\leq n)$.

\smallskip

For a set $H$, $H'$ denotes its \textbf{accumulation points}. Set $H^{(0)}=H$. For $n\in\mathbb{N}$ set $H^{(n)}=\left(H^{(n-1)}\right)'$.
\par $\mathrm{cl}(H)$, $\mathrm{int}(H)$ and $\partial H$ denote the \textbf{closure}, \textbf{interior} and \textbf{boundary} of $H$ respectively.

\smallskip

If $K\subset\mathbb{R}^n$, then $\dim_{\mathrm{H}}(K)$ will denote the \textbf{Hausdorff dimension} of $K$, and $\mu^{d}(K)$ will denote the $d$ dimensional \textbf{Hausdorff measure} of $K$.

\smallskip

$\lambda$ will denote the Lebesgue-measure and also the outer measure as well.

\smallskip

If $I\subset\mathbb{R},\ f:I\to\mathbb{R}$, the \textbf{oscillation} of $f$ on $I$ is
\[\omega_f(I)=\sup\{|f(x)-f(y)|:x,y\in I\},\]
while the oscillation of $f$ at $a\in I$ is
\[\omega_f(a)=\inf\limits_{\delta>0}\omega_f(S(a,\delta))\ \ (=\lim\limits_{\delta\to0+0}\omega_f(S(a,\delta))).\]

\smallskip

If $a\in I\subset\mathbb{R},\ f:I\to\mathbb{R}$, the \textbf{cluster set} $\Lambda_f(a)$ is
\[\Lambda_f(a)=\big\{y\in\mathbb{R}:\exists (x_n)\text{ sequence such that }x_n\in I, x_n\to a\in I, f(x_n)\to y\big\}.\]
Clearly $f(a)\in\Lambda_f(a)$, $\Lambda_f(a)$ is closed and $\omega_f(a)=\mathrm{diam}(\Lambda_f(a))$. 

\smallskip

Let $\mathrm{C}(H)$ denote all continuous functions with domain $H$. 
For ${f,g:H\to\mathbb{R}}$ let us use the following notation: $d_{\sup}(f,g)=\sup\{|f(x)-g(x)|:x\in H\}$.

\section{The bases of the integral}

First, for the integral and related notions, we will need some algebraic structures on $[0,+\infty)\times[-\infty,+\infty]$ and its subsets. We define them in the spirit of the Hausdorff-dimension and measure.

\subsection{Structures on $[0,+\infty)\times[-\infty,+\infty]$ and its subsets}

We will give base definitions on $[0,+\infty)\times[-\infty,+\infty]$, but in several cases we will apply those on subsets using the corresponding subspace structure.

\begin{df}If $(d_1,m_1),(d_2,m_2),\in[0,+\infty)\times[-\infty,+\infty]$, then let
\[(d_1,m_1)\leq(d_2,m_2)\iff d_1<d_2\text{ or }(d_1=d_2\text{ and }m_1\leq m_2),\]
which is the lexicographic order.\qed
\end{df}

\begin{prp}If $K\subset[0,1]\times[-\infty,+\infty]$, then both $\inf K$ and $\sup K$ exist in $[0,1]\times[-\infty,+\infty]$. Similar statement is true in $[0,1]\times[0,+\infty]$.
\end{prp}
\begin{proof}Let $\pi_1$ and $\pi_2$ denote the projection to the first and second coordinate respectively. Let $m=\inf\limits_{k\in K}\pi_1(k),\ M=\sup\limits_{k\in K}\pi_1(k)$. Set
\[i=\begin{cases}
\left(m,+\infty\right)&\text{if }\min\limits_{k\in K}\pi_1(k)\text{ does not exist}\\
\left(m,\ \inf\{\pi_2(k):k\in K,\pi_1(k)=m\}\right)&\text{if }m=\min\limits_{k\in K}\pi_1(k)\text{ exist},
\end{cases}\]
\[s=\begin{cases}
\left(M,-\infty\right)&\text{if }\max\limits_{k\in K}\pi_1(k)\text{ does not exist}\\
\left(M,\ \sup\{\pi_2(k):k\in K,\pi_1(k)=M\}\right)&\text{if }M=\max\limits_{k\in K}\pi_1(k)\text{ exist}.
\end{cases}\]
It can be readily seen that $i=\inf K$ and $s=\sup K$ hold in $[0,1]\times[-\infty,+\infty]$.
\par In $[0,1]\times[0,+\infty]$ only the first line of the definition of $s$ has to be modified to $(M,0)$.
\end{proof}

\begin{df}If $(d_1,m_1),(d_2,m_2),\in[0,+\infty)\times[-\infty,+\infty]$, then let
\[(d_1,m_1)+(d_2,m_2)=\begin{cases}
(d_2,m_2)&\text{if }d_1<d_2\\
(d_1,m_1)&\text{if }d_2<d_1\\
(d_1,m_1+m_2)&\text{if }d_1=d_2\text{ and }m_1+m_2\text{ exists}\\
\text{undefined}&\text{if }d_1=d_2\text{ and }m_1+m_2\text{ is undefined}.\tag*{\qed}
\end{cases}\]
\end{df}

\begin{rem}\label{r1}One can formulate the summation in the following way.
\[(d_1,m_1)+(d_2,m_2)=\left(\max\{d_1,d_2\},\sum_{d_i=\max\{d_1,d_2\}}m_i\right).\tag*{\qed}\]
\end{rem}


In the sequel, we will apply the order topology on $[0,1]\times[0,+\infty]$ i.e. if we do not specify which topology is applied, then we mean the order topology.

\begin{prp}$\langle[0,1]\times[0,+\infty];+\rangle$ is a topological commutative monoid (with $(0,0)$ as identity).
\end{prp}
\begin{proof}Clearly $(d,m)+(0,0)=(d,m)$ and the summation is commutative.
\par We show that $+$ is associative. Let $(d_1,m_1),(d_2,m_2),(d_3,m_3)\in[0,1]\times[0,+\infty]$. Let $d=\max\{d_1,d_2,d_3\}$ and $m=\sum_{d_i=d}m_i$. Evidently both $((d_1,m_1)+(d_2,m_2))+(d_3,m_3)$ and $(d_1,m_1)+((d_2,m_2)+(d_3,m_3))$ equals to $(d,m)$.
\par We show that the summation is continuous. Let $(d,m)=(d_1,m_1)+(d_2,m_2)$. There are two cases. 
\vspace{-0.6pc}\begin{enumerate}\setlength\itemsep{-0.3em}
\item If $d_1=d_2$, then continuity is the consequence of the continuity of the summation on the real numbers.
\item If $d_1<d_2$ and $U$ is a neighborhood of $(d,m)=(d_2,m_2)$, then choose a neighborhood $V$ of $(d_1,m_1)$ and a new neighborhood $U'$ of $(d_2,m_2)$ such that $U'\subset U$ and $V\cap U'=\emptyset$. It can be readily seen that the neighborhoods $V,U'$ will suit i.e. $V+U'=U'\subset U$ holds.\qedhere
\end{enumerate}
\end{proof}

\subsubsection{Sequences and series on $[0,1]\times[-\infty,+\infty]$}

We will need some trivial statements on convergence in $[0,1]\times[-\infty,+\infty]$ later (e.g. for Fatou theorem \ref{tFatou}) that we enumerate here.

\begin{prp}\label{p4}Let $(h_n)$ be a sequence in $[0,1]\times[-\infty,+\infty]$ such that $h_n\leq h_{n+1}$. Then $(h_n)$ is convergent in the order topology, moreover $h_n\to\sup\{h_n:n\in\mathbb{R}\}$.\qed
\end{prp}

\begin{prp}Let $(a_n), (b_n)$ be two sequences in $[0,1]\times[-\infty,+\infty]$ such that $a_n\leq b_n\ (n\in\mathbb{R})$ and $a_n\to a, b_n\to b$. Then $a\leq b$.\qed
\end{prp}

\begin{df}Let $(h_n)$ be a sequence in $[0,1]\times[-\infty,+\infty]$. Set
\[\varliminf\limits_{n\to\infty}h_n=\lim\limits_{n\to\infty}\inf\{h_k:k\geq n\},\ \ \varlimsup\limits_{n\to\infty}h_n=\lim\limits_{n\to\infty}\sup\{h_k:k\geq n\}\tag*{\qed}.\]
\end{df}

\begin{rem}As the sequence $b_k=\inf\{h_k:k\geq n\}$ is increasing, then by \ref{p4} the definition of $\varliminf$ makes sense. Similarly for $\varlimsup$.
\end{rem}

\begin{prp}\label{p5}Let $(a_n)$ be a sequences in $[0,1]\times[-\infty,+\infty]$ such that $a_n\to a$. Then $a=\varliminf\limits_{n\to\infty}a_n=\varlimsup\limits_{n\to\infty}a_n$.\qed
\end{prp}

\begin{prp}\label{p6}Let $(a_n), (b_n)$ be two sequences in $[0,1]\times[-\infty,+\infty]$ such that $a_n\leq b_n\ (n\in\mathbb{R})$ and $a_n\to a$. Then $a\leq\varliminf\limits_{n\to\infty}b_n$.
\end{prp}
\begin{proof}Trivially $\varliminf_{n\to\infty}a_n\leq\varliminf_{n\to\infty}b_n$ and apply \ref{p5}.
\end{proof}

\begin{thm}Let a sequence $((d_i,m_i))$ be given on $[0,1]\times[-\infty,+\infty]$. Let $a_n=\sum_{i=1}^n(d_i,m_i)\ (n\in\mathbb{N})$. Then if $\lim a_n$ exists, then 
$\lim a_n=(d,m)\text{ where }d=\sup\{d_i:i\in\mathbb{N}\},\ m=\sum_{d_i=d}m_i$, i.e.
\[\sum_{i=1}^{\infty}(d_i,m_i)=(d,m)=\left(\sup\{d_i:i\in\mathbb{N}\};\sum_{d_j=\sup\{d_i:i\in\mathbb{N}\}}m_j\right),\]
using the convention that empty sum is $0$.
\end{thm}
\begin{proof}One can simply generalize \ref{r1} using induction and get that 
\[a_n=(d'_n,m'_n)\text{ where }d'_n=\max\{d_1,\dots,d_n\}\text{ and }m'_n=\sum\limits_{\substack{1\leq i\leq n\\d_i=d'_n}}m_i.\]
Clearly $d'_n\to d$. Now there are two cases.
\vspace{-0.6pc}\begin{enumerate}\setlength\itemsep{-0.3em}
\item If for all $n\in\mathbb{N}\ d'_n<d$, then evidently $(d'_n,m'_n)\to(d,0)$.
\item If there is $N\in\mathbb{N}$ such that $d'_N=d$, i.e. $\sup=\max$, then $m'_n=\sum\limits_{1\leq i\leq n, d_i=d}m_i$ for $n\geq N$. In this case $d'_n\to d$ and it can be readily seen that $(d'_n,m'_n)\to(d,m)$.\qedhere
\end{enumerate}
\end{proof}

\subsection{Basic definitions and results}

We define a generalization of (pseudo-)metric.

\begin{df}$d:X\times X\to[0,+\infty)\times[0,+\infty]$ is an \textbf{h-metric} on $X$ if the followings hold.
\vspace{-0.6pc}\begin{enumerate}\setlength\itemsep{-0.3em}
\item $d(x,y)=(0,0)$ iff $x=y$ ($x,y\in X$)
\item $d(x,y)=d(y,x)$ if $x,y\in X$
\item $d(x,z)\leq d(x,y)+d(y,z)$ if $x,y,z\in X$.
\end{enumerate}
We call $\langle X,d\rangle$ an \textbf{h-metric space}.
\newline It is a \textbf{pseudo-h-metric} if in axiom 1, we require $d(x,x)=(0,0)$ only.\qed
\end{df}

\begin{prp}If $d$ is an (pseudo-)h-metric on $X$, then $\pi_1(d(x,y))$ is a pseudo-metric on $X$. If $c\in[0,+\infty]$, then $d'(x,y)=\pi_2(d(x,y))$ is a pseudo-h-metrics on $X'$, where $X'$ defined as a set which satisfies that $x,y\in X'$ implies that $\pi_1(d(x,y))=c$.
\end{prp}
\begin{proof}We just check the triangle inequality in both cases. In the first case, $\pi_1(d(x,y)+d(y,z))=\max\{\pi_1(d(x,y));\pi_1(d(y,z))\}$ gives the statement. In the second case, as the first coordinate is fixed, the second coordinate decides the triangle inequality.
\end{proof}

\begin{ex}If $d$ is an (pseudo-)h-metric on $X$, then $\pi_2(d(x,y))$ is not a pseudo-metric on $X$ in general.
\par Let $d(x,z)=(1,1),\ d(x,y)=(2,0),\ d(y,z)=(2,0)$.\qed
\end{ex}

The next proposition shows that this is really a generalization.

\begin{prp}If $d$ is a (pseudo)-metric on $X$, then $(0,d(x,y))$ is a (pseudo)-h-metric on $X$.\qed
\end{prp}

\begin{df}We define $d^\mathrm{H}:[0,1]\times[0,+\infty]\to[0,1]\times[0,+\infty]$. If $(d_1,m_1),(d_2,m_2),\in[0,1]\times[0,+\infty]$, then set
\[d^\mathrm{H}\big((d_1,m_1),(d_2,m_2)\big)=\begin{cases}
(|d_1-d_2|,0)&\text{if }d_1\ne d_2\\
(0,|m_1-m_2|)&\text{if }d_1=d_2.\tag*{\qed}
\end{cases}\]
\end{df}

\begin{prp}$d^\mathrm{H}$ is a h-metric on $[0,1]\times[0,+\infty]$.
\end{prp}
\begin{proof}Clearly axiom 1 and 2 are satisfied.
\par Let $(d_1,m_1),(d_2,m_2),(d_3,m_3)\in[0,1]\times[0,+\infty]$. To see the triangle axiom, there are 5 cases to be checked. We just enumerate the cases, each of them trivially satisfies the axiom. Cases: $d_1=d_2=d_3$, $d_1\ne d_2\ne d_3\ne d_1$, $d_1=d_2\ne d_3$, $d_1\ne d_2=d_3$, $d_1=d_3\ne d_2$.
\end{proof}

\begin{df}We define $\dm:\mathrm{P}(\mathbb{R})\to\{0\}\times\overline{\mathbb{N}}_0\cup(0,1]\times[0,+\infty]$. If $K\subset\mathbb{R}$, then
\[\dm(K)=(d,m)\text{ where }d=\dimh K,m=\mu^d(K).\tag*{\qed}\]
\end{df}

\begin{prp}\label{p1}If $A,B\subset\mathbb{R}$ and $A\subset B$, then $\dm(A)\leq\dm(B)$.
\end{prp}
\begin{proof}It is because of both dimension and measure are monotone.
\end{proof}

\begin{prp}\label{p2}If $A,B\subset\mathbb{R}$, then $\dm(A\cup B)\leq\dm(A)+\dm(B)$.
\end{prp}
\begin{proof}If $\dimh(A)<\dimh(B)$, then $\dm(A\cup B)=\dm(B)$ and $\dm(A)+\dm(B)=\dm(B)$. If $\dimh(A)=\dimh(B)$, then it is the property of \mbox{(Hausdorff-)} measure.
\end{proof}

\begin{df}\label{dHint}Let $K\subset\mathbb{R}$ be Borel-measurable and $f:K\to\mathbb{R}$ be a Borel-measurable function. Set $D_f=\{x\in K:f(x)\ne 0\}$. Set
\[(\mathrm{H})\int\limits_Kf=(d,m)\in[0,1]\times[-\infty,+\infty],\]
where
\[d=\dimh D_f,\ m=\int\limits_{D_f}f\ d\mu^d,\]
and we use the conventions that $\dimh\emptyset=0$ and $\int_{\emptyset}f\ d\mu^d=0$.
\par It is called the \textbf{Hausdorff-integral} of $f$ on $K$. We call $f$ \textbf{Hausdorff-integrable} if its Hausdorff-integral exists and $m\in\mathbb{R}$.\qed
\end{df}

Clearly the result of the Hausdorff-integral is not a number, it is a pair of numbers (including infinities).

\begin{rem}It is known that the Borel-sets are $\mu^d$-measurable (for any $d$), hence $f$ can be integrated by $\mu^d$. Actually it could have been enough to assume that $f$ is measurable by $\mu^d$ for $d=\dim_{\mathrm{H}}D_f$.\qed
\end{rem}

\begin{rem}If $f,D_f,d,m$ are the ones defined above (\ref{dHint}), then
\[m=\int\limits_{K}f\ d\mu^d,\]
i.e. we could have used that too in the definition.\qed
\end{rem}

\begin{prp}\label{p8}Let $K\subset\mathbb{R}$ be Borel-measurable. If for $f:K\to\mathbb{R}$ the Hausdorff-integral exists on $K$, then it exists on $D_f$ and
\[(\mathrm{H})\int\limits_Kf=(\mathrm{H})\int\limits_{D_f}f\tag*{\qed}.\]
\end{prp}

\begin{rem}If $\lambda(D_f)>0$, then the Hausdorff-integral gives the Lebesgue-integral, more precisely in that case we get that $(\mathrm{H})\int_Kf=(1,m)$ where $m=\int_Kf\ d\lambda$.
\end{rem}

\begin{prp}Let $K\subset\mathbb{R}$ be Borel-measurable. Then 
\[\dm(K)=(\mathrm{H})\int\limits_{\mathbb{R}}\chi_K=(\mathrm{H})\int\limits_K1.\tag*{\qed}\]
\end{prp}

\begin{prp}Let $K\subset\mathbb{R}$ be Borel-measurable and $f:K\to\mathbb{R}$ be Borel-measurable. If $0\leq f$, then its Hausdorff-integral exists on $K$.\qed
\end{prp}

\begin{prp}\label{p7}Let $K\subset\mathbb{R}$ be Borel-measurable and $f:K\to\mathbb{R}$ be a Borel-measurable function. Let $A, B\subset K$ be Borel-measurable and $A\cap B=\emptyset$. Then
\[(\mathrm{H})\int\limits_{A\cup B}f=(\mathrm{H})\int\limits_Af + (\mathrm{H})\int\limits_Bf,\]
whenever all integral exist.
\end{prp}
\begin{proof}For $L\subset K$ let $D^L_f=\{x\in L:f(x)\ne 0\}$ and $d^L=\dimh D^L_f$. Clearly $D^{A\cup B}_f=D^A_f\cup^* D^B_f$ which gives that $d^{A\cup B}=\max\{d^A,d^B\}$. If $d^A\ne d^B$, then we are done. If $d^A=d^B$, then apply the usual measure theory statement.
\end{proof}

We can say more.

\begin{prp}Let $K\subset\mathbb{R}$ be Borel-measurable and $f:K\to\mathbb{R}$ be a Borel-measurable function. Let $(A_n)$ be a sequence of Borel-measurable sets on $K$ such that $A_n\cap A_m=\emptyset\ (n\ne m)$. Let $A=\bigcup_{n=1}^{\infty}A_n$. Then
\[(\mathrm{H})\int\limits_{A}f=\sum_{n=1}^{\infty}(\mathrm{H})\int\limits_{A_n}f,\]
whenever all integral exist.
\end{prp}
\begin{proof}For $L\subset K$ let $D^L_f=\{x\in L:f(x)\ne 0\}$ and $d^L=\dimh D^L_f$. Clearly $D^A_f=\bigcup_{n=1}^{\infty}D^{A_n}_f$ which gives that $d=\sup\{d_n:n\in\mathbb{N}\}$ where $d=d^A, d_n=d^{A_n}$. Let $m_n=\int_{D^{A_n}_f}f\ d\mu^{d_n}$, i.e. $(\mathrm{H})\int\limits_{A_n}f=(d_n,m_n)$. There are two cases.
\vspace{-0.6pc}\begin{enumerate}\setlength\itemsep{-0.3em}
\item If for all $n\in\mathbb{N}\ d_n<d$, then $\mu^d(A_n)=0$ for all $n$, hence $\mu^d(A)=0$ as well, which gives that $(\mathrm{H})\int\limits_{A}f=(d,0)$ and also $\sum_{n=1}^{\infty}(\mathrm{H})\int\limits_{A_n}f=(d,0)$.
\item Suppose now that there is $N\in\mathbb{N}$ such that $d_N=d$. Note that $\mu^d(\bigcup_{d_n<d}A_n)=0$, therefore 
\[(\mathrm{H})\int\limits_{A}f=(\mathrm{H})\int\limits_{\bigcup\limits_{d_n=d}A_n}f=\left(d,\sum_{d_n=d}m_n\right)\]
by the usual measure theory statement. While we also get that
\[\sum_{n=1}^{\infty}(\mathrm{H})\int\limits_{A_n}f=\sum_{n=1}^{\infty}(d_n,m_n)=\left(d,\sum_{d_n=d}m_n\right).\qedhere\]
\end{enumerate} 
\end{proof}

\begin{cor}Let $K\subset\mathbb{R}$ be Borel-measurable and $f:K\to\mathbb{R}$ such that its Hausdorff-integral exists on $K$. Let $f^+=\max\{f,0\}, f^-=\min\{f,0\}$. Then
\[(\mathrm{H})\int\limits_Kf=(\mathrm{H})\int\limits_Kf^+  +  (\mathrm{H})\int\limits_Kf^-.\]
\end{cor}
\begin{proof}Obviously $D_f=D_{f^+}\cup^*D_{f^-}$. Then
\[(\mathrm{H})\int\limits_Kf=(\mathrm{H})\int\limits_{D_f}f=(\mathrm{H})\int\limits_{D_{f^+}\cup D_{f^-}}f=(\mathrm{H})\int\limits_{D_{f^+}}f+(\mathrm{H})\int\limits_{D_{f^-}}f=\]
\[(\mathrm{H})\int\limits_{D_{f^+}}f^++(\mathrm{H})\int\limits_{D_{f^-}}f^-=(\mathrm{H})\int\limits_Kf^++(\mathrm{H})\int\limits_Kf^-\]
by propositions \ref{p8} and \ref{p7}.
\end{proof}

\begin{prp}\label{p11}Let $K\subset\mathbb{R}$ be Borel-measurable. If $f:K\to\mathbb{R}$ is Hausdorff-integrable, $(\mathrm{H})\int_Kf=(d,m)$ and $c\in\mathbb{R},\ c\ne 0$, then $cf$ is also Hausdorff-integrable and $(\mathrm{H})\int_Kcf=(d,cm)$.\qed
\end{prp}

\begin{thm}\label{p10}Let $K\subset\mathbb{R}$ be Borel-measurable. If $f,g:K\to\mathbb{R}$ are Hausdorff-integrable such that $0\leq f,g$, then $f+g$ is also Hausdorff-integrable and
\[(\mathrm{H})\int\limits_Kf+g=(\mathrm{H})\int\limits_Kf  +  (\mathrm{H})\int\limits_Kg.\]
\end{thm}
\begin{proof}Evidently $D_{f+g}=D_f\cup D_g$. There are 3 cases.
\vspace{-0.6pc}\begin{enumerate}\setlength\itemsep{-0.3em}
\item If $\dimh D_f=\dimh D_g$, then apply the usual measure theory statement.
\item If $\dimh D_g<\dimh D_f$, then by propositions \ref{p8}
\[(\mathrm{H})\int\limits_Kf+g=(\mathrm{H})\int\limits_{D_f\cup D_g}f+g=(\mathrm{H})\int\limits_{D_f-D_g}f+g=(\mathrm{H})\int\limits_{D_f}f=(\mathrm{H})\int\limits_Kf  +  (\mathrm{H})\int\limits_Kg\]
because $\mu^d(D_g)=0$ for $d=\dimh D_f=\dimh D_f\cup D_g$.
\item If $\dimh D_f<\dimh D_g$, then similar argument works.\qedhere
\end{enumerate}
\end{proof}

\begin{ex}The previous theorem does not hold in general.
\par Let $K=[0,1]$ and $f(x)\equiv 1,\ g(x)\equiv -1$. Then $(\mathrm{H})\int\limits_Kf+g=(\mathrm{H})\int\limits_K0=(0,0)$ while $(\mathrm{H})\int\limits_Kf=(1,1)$ and $(\mathrm{H})\int\limits_Kg=(1,-1)$.\qed
\end{ex}

\begin{thm}\label{p9}Let $K\subset\mathbb{R}$ be Borel-measurable. If $f,g:K\to\mathbb{R}$ are Hausdorff-integrable and $0\leq f\leq g$, then
\[(\mathrm{H})\int\limits_Kf\leq(\mathrm{H})\int\limits_Kg.\]
\end{thm}
\begin{proof}First note that $D_f\subset D_g$. There are two cases.
\vspace{-0.6pc}\begin{enumerate}\setlength\itemsep{-0.3em}
\item If $\dimh D_f<\dimh D_g$, then we are done.
\item If $\dimh D_f=\dimh D_g$, then it is reduced to the usual measure theory statement.\qedhere
\end{enumerate}
\end{proof}

\begin{ex}The previous theorem does not hold in general.
\par Let $K=[0,1]$ and $f(x)\equiv -1,\ g(x)=0$ for $x\ne 0$ and $g(0)=1$. Then $(\mathrm{H})\int\limits_Kf=(1,-1)$ while $(\mathrm{H})\int\limits_Kg=(0,1)$.\qed
\end{ex}

\begin{cor}Let $K\subset L\subset\mathbb{R}$ be Borel-measurable sets. If ${f:L\to\mathbb{R}}$ Borel-measurable and $0\leq f$, then
\[(\mathrm{H})\int\limits_Kf\leq(\mathrm{H})\int\limits_Lf.\]
\end{cor}
\begin{proof}Apply theorem \ref{p9} for functions $f\chi_K$ and $f$.
\end{proof}

\begin{prp}Let $K\subset\mathbb{R}$ be Borel-measurable and $f:K\to\mathbb{R}$ be Borel-measurable function such that $(\mathrm{H})\int\limits_Kf=(0,m)$ and $m\in\mathbb{R}$. Then $D_f$ is countable and 
\[m=\sum_{x\in D_f}f(x)\]
where the result of the sum is independent of the order of the terms.
\end{prp}
\begin{proof}If the integral exists and its first coordinate equals to $0$, then $\dimh D_f=0$. If $D_f$ is uncountable, then there is $\varepsilon>0$ such that $\{x\in K:|f(x)|>\varepsilon\}$ is uncountable which gives that either $\int_Kf\ d\mu^0$ equals to $\pm\infty$ or the integral does not exist. Therefore $D_f$ has to be countable.
\par As $\mu^0$ is the counting measure and $\mu^0(\{x\})=1$ for $x\in K$, we get that
\[\int\limits_Kf\ d\mu^0=\int\limits_{D_f}f\ d\mu^0=\sum\limits_{x\in D_f}\int\limits_{\{x\}}f\ d\mu^0=\sum\limits_{x\in D_f}f(x).\qedhere\]
\end{proof}

\subsection{Integrals of series}

We show that the analog of the Beppo-Levi and Fatou theorems hold here, but only for non-negative functions.

\begin{thm}\label{tBeppoLevi}Let $K\subset\mathbb{R}$ be Borel-measurable and $(f_n)$ be a sequence of non-negative Borel-measurable functions on $K$ such that $f_n\leq f_{n+1}\ (n\in\mathbb{N})$ and $f_n\to f$. Then
\[(\mathrm{H})\int\limits_Kf_n\to(\mathrm{H})\int\limits_Kf\]
in the order topology on $[0,1]\times[0,+\infty]$.
\end{thm}
\begin{proof}Let 
\[(\mathrm{H})\int\limits_Kf_n=(d_n,m_n),\ \ \ (\mathrm{H})\int\limits_Kf=(d,m).\]
Let $D_n=D_{f_n}, D=D_f$. Then $d_n=\dimh D_n, d=\dimh D$. 
\par By $f_n\leq f_{n+1}$ we get that $D_n\subset D_{n+1}$ for all $n\in\mathbb{N}$. Furthermore by $f_n\to f$, we have $D=\bigcup_{n=1}^{\infty}D_n$. Those give that $d_n\leq d_{n+1}$ and $d_n\to d$.
\par There are two cases.
\vspace{-0.6pc}\begin{enumerate}\setlength\itemsep{-0.3em}
\item For all $n\in\mathbb{N}$, $d_n<d$ holds. But then $\mu^d(D_n)=0$ which gives that $\mu^d(D)=0$ which means that $m=0$. Hence $(d_n,m_n)\to(d,0)$. 
\item There is $N\in\mathbb{N}$ such that $n>N$ implies that $d_n=d$. Then the usual Beppo-Levi theorem yields that $(d_n,m_n)\to(d,m)$, because for $n>N$
\[\int\limits_{D_n}f_n\ d\mu^{d_n}=\int\limits_Kf_n\ d\mu^d\to\int\limits_Kf\ d\mu^d=\int\limits_Df\ d\mu^d.\qedhere\]
\end{enumerate}
\end{proof}

\begin{ex}\label{e1}If we allow negative functions as well, then the Beppo-Levi theorem does not hold in general.
\par Let $f_n:\mathbb{R}_0^+\to\mathbb{R}, f:\mathbb{R}_0^+\to\mathbb{R}$,
\[f_n(x)=\begin{cases}
0&\text{if }\frac{1}{n}<x\\
-1&\text{if }0\leq x\leq \frac{1}{n},
\end{cases}\]
\[f(x)=\begin{cases}
0&\text{if }x>0\\
-1&\text{if }x=0.
\end{cases}\]
Obviously $f_n\leq f_{n+1}\ (n\in\mathbb{N})$ and $f_n\to f$, and $(\mathrm{H})\int\limits_Kf_n=(1,-\frac{1}{n}),\ \ (\mathrm{H})\int\limits_Kf=(0,-1)$.\qed
\end{ex}

\begin{ex}\label{e2}The Beppo-Levi theorem does not hold for decreasing non-negative functions in general.
\par Simply use functions $-f_n\ (n\in\mathbb{N})$ and $-f$ from \ref{e1}.\qed
\end{ex}

\begin{thm}\label{tFatou}Let $K\subset\mathbb{R}$ be Borel-measurable and $(f_n)$ be a sequence of non-negative Borel-measurable functions on $K$ such that $f_n\to f$. Then
\[(\mathrm{H})\int\limits_Kf\leq\varliminf\limits_{n\to\infty}(\mathrm{H})\int\limits_Kf_n.\]
\end{thm}
\begin{proof}Let $g_n=\inf\{f_k:k\geq n\}$. Then $0\leq g_n,\ g_n\leq g_{n+1}\ (n\in\mathbb{N})$ and $g_n\to f$. Then by theorem \ref{tBeppoLevi}
\[(\mathrm{H})\int\limits_Kg_n\to(\mathrm{H})\int\limits_Kf.\]
Clearly $g_n\leq f_n$ which yields that
\[(\mathrm{H})\int\limits_Kg_n\leq(\mathrm{H})\int\limits_Kf_n.\]
by proposition \ref{p9}. By proposition \ref{p6} those two statements together give the claim.
\end{proof}

\begin{rem}The analog of Lebesgue-theorem does not hold in general as we showed it in examples \ref{e1} and \ref{e2}. Similarly the Fatou theorem does not hold if we allow negative functions as well.
\end{rem}

\section{Applications}

Now we define an h-metric on $\mathrm{P}(\mathbb{R})$ which is a refinement of the usual pseudo-metric on $\mathrm{P}(\mathbb{R})$, namely $\lambda(A\triangle B)$. 

\begin{df}If $A,B\subset\mathbb{R}$, then set
\[d_s(A,B)=\dm(A\triangle B).\tag*{\qed}\]
\end{df}

\begin{prp}$d_s$ is an h-metric on $\mathrm{P}(\mathbb{R})$.
\end{prp}
\begin{proof}Clearly $d_s(A,A)=(0,0)$ and $d_s(A,B)=d_s(B,A)$.
\par If $d_s(A,B)=(0,0)$, then $\dimh A\triangle B=0$ and $\mu^0(A\triangle B)=0$, but as $\mu^0$ is the counting measure, we get that $A=B$.
\par To see the triangle inequality, note that $A\triangle C\subset(A\triangle B)\cup(B\triangle C)$ holds, and then apply propositions \ref{p1} and \ref{p2}.
\end{proof}

Now we define an h-metric on Borel-measurable functions which is a refinement of the usual pseudo-metric. 

\begin{df}Let $K\subset\mathbb{R}$ be Borel-measurable. Then set
\[{\cal L}_{\mathrm{H}}(K)=\{f:K\to\mathbb{R}\ :\ |f|\text{ is Hausdorff-integrable}\}.\tag*{\qed}\]
\end{df}

\begin{prp}Let $K\subset\mathbb{R}$ be Borel-measurable. If $c\in\mathbb{R},\ f,g\in{\cal L}_{\mathrm{H}}(K)$, then $cf\in{\cal L}_{\mathrm{H}}(K)$ and $f+g\in{\cal L}_{\mathrm{H}}(K)$.
\end{prp}
\begin{proof}The first statement is a consequence of \ref{p11}.
\par To see the second statement, we get that
\[(\mathrm{H})\int\limits_K|f+g|\leq(\mathrm{H})\int\limits_K|f|+|g|=(\mathrm{H})\int\limits_K|f|+(\mathrm{H})\int\limits_K|g|\]
by \ref{p9} and \ref{p10}.
\end{proof}

\begin{df}Let $K\subset\mathbb{R}$ be Borel-measurable and $f,g\in{\cal L}_{\mathrm{H}}(K)$. Set
\[d_{\mathrm{H}}(f,g)=d(f,g)=(\mathrm{H})\int\limits_K|f-g|.\]
\end{df}

\begin{thm}Let $K\subset\mathbb{R}$. Then $d(f,g)$ is an h-metric on ${\cal L}_{\mathrm{H}}(K)$.
\end{thm}
\begin{proof}Evidently $d(f,f)=(0,0)$ and $d(f,g)=d(g,f)$.
\par If $d(f,g)=(0,0)$, then $\dimh D_{|f-g|}=0$ and $\int\limits_{D_{|f-g|}}|f-g|\ d\mu^0=0$. But as $\mu^0$ is the counting measure, it gives that $f=g$ everywhere.
\par To prove the triangle inequality, let $f,g,h\in{\cal L}_{\mathrm{H}}(K)$. We have to show that
\[(\mathrm{H})\int\limits_K|f-h|\leq(\mathrm{H})\int\limits_K|f-g|+(\mathrm{H})\int\limits_K|g-h|.\]
Let 
\[d_1=\dimh D_{|f-h|},\ d_2=\dimh D_{|f-g|},\ d_3=\dimh D_{|g-h|}.\]
If $|f(x)-h(x)|\ne 0$, then either $|f(x)-g(x)|\ne 0$ or $|g(x)-h(x)|\ne 0$, which gives that $D_{|f-h|}\subset D_{|f-g|}\cup D_{|g-h|}$. 
We get that either $d_1\leq d_2$ or $d_1\leq d_3$. Then we have the following cases.
\vspace{-0.6pc}\begin{enumerate}\setlength\itemsep{-0.3em}
\item If either $d_1<d_2$ or $d_1<d_3$, then we are done.
\item If $d_1=d_2=d_3$, then the statement is reduced to the usual measure theory statement. 
\item If $d_3<d_1=d_2$, then $\dimh D_{|g-h|}<d_1$, hence $\mu^{d_1}(D_{|g-h|})=0$ which gives that $\mu^{d_1}(D_{|f-h|})=\mu^{d_1}(D_{|f-h|}-D_{|g-h|})$. But 
\[D_{|f-h|}-D_{|g-h|}=D_{|f-h|}\cap(K-D_{|g-h|})=D_{|f-g|}\]
because $K-D_{|g-h|)}=\{x\in K:g(x)=h(x)\}$. Therefore we get that
\[\int\limits_{D_{|f-h|}}|f-h|\ d\mu^{d_1}=\int\limits_{D_{|f-h|}-D_{|g-h|}}|f-h|\ d\mu^{d_1}=\]
\[\int\limits_{D_{|f-g|}}|f-g|\ d\mu^{d_1}=\int\limits_{D_{|f-g|}}|f-g|\ d\mu^{d_2},\]
hence
\[(\mathrm{H})\int\limits_K|f-h|=(\mathrm{H})\int\limits_K|f-g|\leq(\mathrm{H})\int\limits_K|f-g|+(\mathrm{H})\int\limits_K|g-h|.\]
\item If $d_2<d_1=d_3$, then a similar argument to the previous one can work.\qedhere
\end{enumerate}
\end{proof}

\begin{rem}We also gained that while $L_1$ was only a pseudo-metric space, ${\cal L}_{\mathrm{H}}(K)$ is now an h-metric space (not just pseudo).
\end{rem}

\begin{df}Let $\langle X,d\rangle$ be an h-metric space. The \textbf{topology} induced by $d$ is the following. The neighborhood of a point $x\in X$ is $\{S(x,(\varepsilon_1,\varepsilon_2)):(\varepsilon_1,\varepsilon_2)>(0,0)\}$, where $S(x,(\varepsilon_1,\varepsilon_2))=\{y\in X:d(x,y)<(\varepsilon_1,\varepsilon_2)\}$.\qed
\end{df}

\begin{prp}Let $\langle X,d\rangle$ be an h-metric space. In the topology induced by $d$, the neighborhood of a point $x\in X$ is $\{S(x,\varepsilon):\varepsilon>0\}$, where $S(x,\varepsilon)=\{y\in X:d(x,y)<(0,\varepsilon)\}$.
\par In the topology induced by $d$, $x_n\to x$ iff $d(x_n,x)\to(0,0)$.\qed
\end{prp}

\begin{df}Let $\langle X,d\rangle$ be an h-metric space. We call a sequence $(x_n)$ \textbf{Cauchy-sequence}, if for all $(\varepsilon_1,\varepsilon_2)>(0,0)$ there is $N\in\mathbb{N}$ such that $n,m>N$ implies that $d(x_n,x_m)<(\varepsilon_1,\varepsilon_2)$.
\end{df}

\begin{prp}\label{p12}Let $\langle X,d\rangle$ be an h-metric space. Then a sequence $(x_n)$ is Cauchy iff for all $\varepsilon>0$ there is $N\in\mathbb{N}$ such that $n,m>N$ implies that $d(x_n,x_m)<(0,\varepsilon)$.\qed
\end{prp}

\begin{df}Let $\langle X,d\rangle$ be an h-metric space. We call it \textbf{complete} if every Cauchy-sequence is convergent.
\end{df}

Now we prove the analog of Riesz-Fischer theorem.

\begin{lem}\label{l1}Let $K\subset\mathbb{R}$ be Borel-measurable, $f:K\to\mathbb{R}$ be a Borel-measurable function and $0\leq f$. If $\int\limits_Kf\ d\mu^0$ is finite, then $\dimh D_{f}=0$.
\end{lem}
\begin{proof}We show that if $\int\limits_Kf\ d\mu^0$ is finite, then $D_{f}$ is countable, which immediately gives that $\dimh D_{f}=0$. Suppose the contrary that $D_{f}$ is uncountable. Then there is $\varepsilon>0$ such that $\{x\in K:f(x)>\varepsilon\}$ is uncountable too. This would yield that $\int\limits_Kf\ d\mu^0=+\infty$ which is a contradiction.
\end{proof}

\begin{thm}Let $K\subset\mathbb{R}$ be Borel-measurable. Then the h-metric space $\langle{\cal L}_{\mathrm{H}}(K), d\rangle$ is complete.
\end{thm}
\begin{proof}Let $(f_n)$ be a Cauchy-sequence in ${\cal L}_{\mathrm{H}}(K)$. Proposition \ref{p12} gives that for all $\varepsilon>0$ there is $N\in\mathbb{N}$ such that $n,m>N$ implies that $\dimh D_{|f_n-f_m|}=0$ and $\int\limits_K|f_n-f_m|\ d\mu^0<\varepsilon$. Clearly this means that $(f_n)$ is a Cauchy-sequence in $L_1(K,\mu^0)$. By the Riesz-Fischer theorem, there is an $f\in L_1(K,\mu^0)$ such that $f_n\to f$ in $L_1(K,\mu^0)$ which is
\begin{equation}\label{eq1}\int\limits_K|f-f_n|\ d\mu^0\to 0.\end{equation}
\par Lemma \ref{l1} and $f\in L_1(K,\mu^0)$ give that $\dimh D_{|f|}=0$, which altogether implies that $f\in{\cal L}_{\mathrm{H}}(K)$.
\par By (\ref{eq1}), there is $N\in\mathbb{N}$ such that $n>N$ implies that $\int\limits_K|f-f_n|\ d\mu^0<1$, and by lemma \ref{l1}, $\dimh D_{|f-f_n|}=0$. These altogether gives that 
\[(\mathrm{H})\int\limits_K|f-f_n|\to(0,0)\]
i.e. $f_n\to f$ in ${\cal L}_{\mathrm{H}}(K)$.
\end{proof}

\begin{rem}We remark that one can also provide a similar proof for the previous theorem to the usual proof of the Riesz-Fischer theorem, namely following the same steps and applying the corresponding versions of the Beppo-Levi and Fatou theorems, noting that only those versions needed that regards for positive functions, and earlier we provided those (theorems \ref{tBeppoLevi} and \ref{tFatou}).
\end{rem}

\subsection{Applications for deficiency measurement}

We are now going to provide five applications of the previously built theory, three for Hausdorff-integral, and two for set measurement. All of them related to deficiency measurement (see \cite{lamd}), which is we want to measure how much a certain property is not fulfilled for some objects. Here we will just define the deficiency measurements with a little explanation, and we will not add reasonable statements on them as our purpose is not that currently.

\smallskip

The notation that we introduced for a deficiency measurement is the following. We have an object $O$, a property $P$ (that is relevant for the object) and a measurement type $t$, and we want to measure how much object $O$ does not satisfy property $P$ when applying measurement type $t$. The result is denoted by $\boldsymbol{\mathrm{defi}(O,P,t)}$ that is often a number in $\overline{\mathrm{R}}_0^+$, but can also be a value in an ordered set with minimum element such that if the property is fulfilled for the object, then the minimum value is taken.

\medskip

First we deal with the overall continuity of a real function. Let $f:\mathbb{R}\to\mathbb{R}$ be a function. As $\omega_f(x)$ naturally measures how much $f$ is not continuous at $x\in\mathbb{R}$, one way of measurement is to ''sum'' those values by an integral. However usually the only problem is that our function is discontinuous on only a set with Lebesgue-measure $0$, hence using the Lebesgue-integral we got $0$ and cannot distinguish it from a continuous function. But if we apply the Hausdorff-integral, we get a much greater granularity.
\[\mathrm{defi}(f,\text{''continuity''},\text{''h-int osc''})=(\mathrm{H})\int\limits_{\mathbb{R}}\omega_f.\]

Another option is to measure the distance of $f$ from the set of all continuous functions. Again if we used the distance in e.g. $L_1$, then many functions would be measured to ''continuous type''. Using the distance in ${\cal L}_{\mathrm{H}}(K)$, which is defined by Hausdorff-integral, we get a much usable result.
\[\mathrm{defi}(f,\text{''continuity''},\text{''dist by ''}d_{\mathrm{H}})=\inf\big\{d_{\mathrm{H}}(f,g) : g\in \mathrm{C}(\mathbb{R})\big\},\]

One more option is to use the sizes of the cluster sets, and consider their ''maximum'' size as a measurement. Using our new set measurement $\dm(K)$, we get a much detailed result than using e.g. Lebesgue-measure.
\[\mathrm{defi}(f,\text{''continuity''},\text{''sup dm cluster''})=\sup\left\{\dm(\Lambda_f(x)):x\in\mathbb{R}\right\}.\]

\medskip

For measuring how much a function $f:\mathbb{R}\to\mathbb{R}$ is not even, we can reflect its graph for the $y$-axis and compare those two function somehow. Again comparing them in $L_1$ is one option, but the Hausdorff-integral can be more handy here too.
\[\mathrm{defi}(f,\text{''even''},\text{''h-int reflect''})=(\mathrm{H})\int\limits_{\mathbb{R}}|f(x)-f(-x)|.\]

\medskip

Let $H\subset\mathbb{R}^2$. We want to measure how much the set $H$ is not convex. A natural option can be to compare its convex hull $\mathrm{Conv}(H)$ to the set itself, and some measurement on the difference of two sets can work. Again, here we apply our new tool for set measurement.
\[\mathrm{defi}(H,\text{''convex''},\text{''dm diff''})=\dm(\mathrm{Conv}(H)-H).\]


\medskip

{\footnotesize

\noindent

\noindent email: alosonczi1@gmail.com\\
}
\end{document}